\documentclass[a4paper,11pt,fleqn]{article}

\usepackage{amsmath,epsfig}

\usepackage{amssymb}
\usepackage{enumerate}
\usepackage{theorem}
\usepackage{pifont}   
\usepackage{euscript}
\usepackage{pstricks}
\usepackage{xcolor}

\title{A forward-backward view of some primal-dual optimization
methods in image recovery}

\author{P. L. Combettes,$\!^1$ L. Condat,$\!^2$ J.-C. Pesquet,$\!^3$
and B. C. V\~u\,$^4$ \footnote{This work was supported by the CNRS
MASTODONS project (grant 2013 MesureHD).}
\\[5mm]
\small
\small $\!^1$Sorbonne Universit\'es -- UPMC Univ. Paris 06\\
\small Laboratoire Jacques-Louis Lions\\
\small Paris, France\\
\small \ttfamily{plc@math.jussieu.fr}
\\[5mm]
\small
\small $\!^3$University of Grenoble--Alpes, GIPSA-lab
\small GIPSA-lab\\
\small St Martin d'H\`eres, France\\
\small \ttfamily{Laurent.Condat@gipsa-lab.grenoble-inp.fr}
\\[5mm]
\small
\small $\!^3$Universit\'e Paris-Est\\
\small Laboratoire d'Informatique Gaspard Monge -- CNRS UMR 8049\\
\small Champs sur Marne, France\\
\small \ttfamily{jean-christophe.pesquet@univ-paris-est.fr}
\\[5mm]
\small
\small $\!^3$LCSL\\
\small Istituto Italiano di Tecnologia and MIT\\
\small Genova, Italy\\
\small \ttfamily{Cong.Bang@iit.it}
}

\date{~}

\renewcommand{\leq}{\ensuremath{\leqslant}}
\renewcommand{\geq}{\ensuremath{\geqslant}}
\renewcommand{\le}{\ensuremath{\leqslant}}

\newcommand{\minimize}[2]{\ensuremath{\underset{\substack{{#1}}}%
{\mathrm{minimize}}\;\;#2 }}

\newcommand{\menge}[2]{\big\{{#1} \mid {#2}\big\}}

\newcommand{\emp}{\ensuremath{{\varnothing}}}

\newcommand{\scal}[2]{\left\langle{#1}\mid {#2} \right\rangle}

\newcommand{\infconv}{\ensuremath{\mbox{\footnotesize$\,\square\,$}}}

\newcommand{\HH}{\ensuremath{\mathcal H}}
\newcommand{\GG}{\ensuremath{\mathcal G}}
\newcommand{\BL}{\ensuremath{\EuScript B}\,}
\newcommand{\SL}{\ensuremath{\EuScript S}\,}
\newcommand{\BP}{\ensuremath{\EuScript P}}

\newcommand{\KKK}{\ensuremath{\boldsymbol{\mathcal K}}}
\newcommand{\GGG}{\ensuremath{\boldsymbol{\mathcal G}}}

\newcommand{\RR}{\ensuremath{\mathbb R}}

\newcommand{\RP}{\ensuremath{\left[0,+\infty\right[}}
\newcommand{\RPP}{\ensuremath{\,\left]0,+\infty\right[}}
\newcommand{\RX}{\ensuremath{\,\left]-\infty,+\infty\right]}}

\newcommand{\RXX}{\ensuremath{\left[-\infty,+\infty\right]}}
\newcommand{\NN}{\ensuremath{\mathbb N}}
\newcommand{\dom}{\ensuremath{\operatorname{dom}}}

\newcommand{\prox}{\ensuremath{\operatorname{prox}}}

\newcommand{\argmin}{\ensuremath{\operatorname{argmin}}}
\newcommand{\ran}{\ensuremath{\operatorname{ran}}}
\newcommand{\zer}{\ensuremath{\operatorname{zer}}}
\newcommand{\gra}{\ensuremath{\operatorname{gra}}}

\newcommand{\Id}{\ensuremath{\operatorname{Id}}}

\newcommand{\weakly}{\ensuremath{\rightharpoonup}}

\newcommand{\pinf}{\ensuremath{+\infty}}
\newtheorem{theorem}{Theorem}[section]

\newtheorem{proposition}[theorem]{Proposition}

\theoremstyle{plain}{\theorembodyfont{\rmfamily}
}
\theoremstyle{plain}{\theorembodyfont{\rmfamily}
}
\theoremstyle{plain}{\theorembodyfont{\rmfamily}
}
\theoremstyle{plain}{\theorembodyfont{\rmfamily}
}
\theoremstyle{plain}{\theorembodyfont{\rmfamily}
\newtheorem{problem}[theorem]{Problem}}
\theoremstyle{plain}{\theorembodyfont{\rmfamily}
}
\theoremstyle{plain}{\theorembodyfont{\rmfamily}
}

\begin{document}

\maketitle

\vspace*{-1cm}
\begin{abstract}
A wide array of image recovery problems can be abstracted into the
problem of minimizing a sum of composite convex functions in a
Hilbert space. To solve such problems, primal-dual proximal
approaches have been developed which provide efficient solutions to
large-scale optimization problems. The
objective of this paper is to show that a number of existing
algorithms can be derived from a general form of the
forward-backward algorithm applied in a suitable product space.
Our approach also allows us to develop useful extensions of
existing algorithms by introducing a variable metric. An
illustration to image restoration is provided.
\end{abstract}


\section{Introduction}
Many image recovery problems can be formulated 
in Hilbert spaces $\HH$ and $(\GG_i)_{1\le i \le m}$
as structured optimization problems of the form
\begin{equation}
\label{e:prob1}
\underset{x\in\HH}{\mathrm{minimize}}\;\;\sum_{i=1}^m g_i(L_i x),
\end{equation}
where, for every $i\in \{1,\ldots,m\}$, $g_i$ is a proper lower
semicontinuous convex function from $\GG_i$ to $\RX$ and $L_i$ is a
bounded linear operator from $\HH$ to $\GG_i$. For example, the
functions $(g_i \circ L_i)_{1\leq i \leq m}$ may model
data fidelity terms, smooth or nonsmooth measures of regularity, or hard
constraints on the solution. In recent years,
many algorithms have been developed to solve such a problem by
taking advantage of recent advances in convex optimization,
especially in the development of proximal tools (see
\cite{Combettes_P_2010_inbook_proximal_smsp,Sra_S_2012_book_optimization_ml}
and the references therein). 
In image processing, however,
solving such a problem still poses a number of conceptual and
numerical challenges. First of all, one often looks for methods
which have the ability to split the problem by activating each of
the functions through elementary processing steps which can be
computed in parallel. This makes it possible to reduce the complexity
of the original problem and to benefit from existing parallel
computing architectures.  Secondly, it is often useful to design
algorithms which can exploit, in a flexible manner, the 
structure of the problem. In particular,
some of the functions may be Lipschitz differentiable in which case
they should be exploited through their gradient rather than 
through their proximity operator, which is usually harder to 
implement (examples of proximity operators with closed-form
expression can be found in 
\cite{Chaux_C_2007_j-ip_variational_ffbip,%
Combettes_P_2010_inbook_proximal_smsp}). 
In some problems, the functions $(g_i)_{1\leq i\leq m}$ can be 
expressed as the infimal convolution of simpler functions
(see \cite{Combettes_P_2013_siam-opt_sys_smi}
and the references therein).
Last but not least, in image recovery, the operators
$(L_i)_{1\leq i \leq m}$ may be of very large size so that their
inversions are costly (e.g., in reconstruction problems). Finding
algorithms which do not require to perform inversions of these
operators is thus of paramount importance. 

Note that all the existing convex optimization algorithms do not
have these desirable properties. For example, the Alternating
Direction Method of Multipliers (ADMM)
\cite{Fortin_M_1983_book_augmented_lmansbvp,Figueiredo_M_2009_ssp_Deconvolution_opiuvsaalo,Goldstein_T_2009_j-siam-is_split_bml}
requires a stringent assumption of invertibility of the involved
linear operator.
Parallel versions of ADMM
\cite{Setzer_S_2009_j-jvcir_deblurring_pibsbt} and related Parallel
Proximal Algorithm (PPXA)
\cite{Combettes_PL_2008_j-ip_proximal_apdmfscvip,Pesquet_J_2012_j-pjpjoo_par_ipo}
usually necessitate a linear inversion to be performed at each
iteration. Also, early primal-dual algorithms 
\cite{Briceno_L_2011_j-siam-opt_mon_ssm,Chambolle_A_2010_first_opdacpai,Chen_G_1994_j-mp_pro_bdm,Combettes_P_2010_j-svva_dualization_srp,Esser_E_2010_j-siam-is_gen_fcf,He_B_2012_j-siam-is_conv_apd}
did not make it possible to handle smooth functions through their gradients.
Only recently, have primal-dual methods been proposed with this
feature. Such work was initiated in
\cite{Combettes_P_2012_j-svva_pri_dsa} in the line of
\cite{Briceno_L_2011_j-siam-opt_mon_ssm} and subsequent developments
can be found in
\cite{Bec_S_2014_j-nonlinear-conv-anal_alg_sps,Bot_R_2013_jmiv_conv_primal_apd,Chen_P_2013_j-inv-prob_prim_dfp,Combettes_P_2013_siam-opt_sys_smi,Condat_L_2013_j-ota-primal-dsm,Repetti_A_2012_p-eusipco_penalized_wlsardcsdn,Vu_B_2013_j-acm_spl_adm}.
As will be seen in the present paper, another advantage of these
approaches is that they can be coupled with variable metric
strategies which can potentially accelerate their convergence.

In Section~\ref{sec:2}, we provide some background on convex
analysis and monotone operator theory. In Section~\ref{sec:3}, we
introduce a general form of the forward-backward algorithm which
uses a variable metric. This algorithm is 
employed in Section~\ref{sec:4} to
develop a versatile family of primal-dual proximal methods. Several
particular instances of this framework are discussed. Finally,
we provide illustrating numerical results in Section~\ref{sec:5}.

\section{Notation and background}
\label{sec:2}
Monotone operator theory \cite{Bauschke_H_2011_book_con_amo}
provides a both insightful and elegant framework for dealing with
convex optimization problems and developing new solution algorithms
that could not be devised using purely variational tools. 
We summarize a number of related concepts that will be needed.

Throughout, $\HH$, $\GG$, and $(\GG_i)_{1\leq i\leq m}$ are real 
Hilbert spaces. We denote the scalar product of a Hilbert space by 
$\scal{\cdot}{\cdot}$ and the associated norm by $\|\cdot\|$. 
The symbol $\weakly$ 
denotes weak convergence,\footnote{In a finite dimensional space, weak convergence is equivalent to strong convergence.} 
and $\Id$ denotes the identity operator.
We denote by $\BL(\HH,\GG)$ the space of bounded linear operators 
from $\HH$ to $\GG$, we set $\SL(\HH)=\menge{L\in\BL(\HH,\HH)}{L=L^*}$, where $L^*$ denotes the
adjoint of $L$. The Loewner partial ordering on $\SL(\HH)$ is 
denoted by $\succcurlyeq$.
For every $\alpha\in\RP$, we set
$
\BP_{\alpha}(\HH)=\menge{U\in\SL(\HH)}{U\succcurlyeq\alpha\Id},
$
and we denote by $\sqrt{U}$ the square root of 
$U\in\BP_{\alpha}(\HH)$. Moreover, for every 
$U\in\BP_\alpha(\HH)$ and $\alpha>0$, we define the norm 
$\|x\|_U=\sqrt{\scal{Ux}{x}}$.

We denote by $\GGG=\GG_1\oplus\cdots\oplus\GG_m$ the Hilbert direct 
sum of the Hilbert spaces $(\GG_i)_{1\leq i\leq m}$, i.e., their 
product space equipped with the scalar product 
$\colon(\boldsymbol{x},\boldsymbol{y})\mapsto
\sum_{i=1}^m\scal{x_i}{y_i}$
where ${\boldsymbol x}=(x_i)_{1\leq i\leq m}$ and 
${\boldsymbol y}=(y_i)_{1\leq i\leq m}$ denote generic elements in
$\GGG$.

Let $A\colon\HH\to 2^{\HH}$ be a set-valued operator.
We denote by $\gra A=\menge{(x,u)\in\HH\times\HH}{u\in Ax}$
the graph of $A$, by
$\zer A=\menge{x\in\HH}{0\in Ax}$ the set of zeros 
of $A$, and by  $\ran A=\menge{u\in\HH}{(\exists\; x\in\HH)\;u\in Ax}$ 
its range. 
The inverse of $A$ is
$A^{-1}\colon\HH\mapsto 2^{\HH}\colon u\mapsto 
\menge{x\in\HH}{u\in Ax}$, and the resolvent of $A$ is
$J_A=(\Id+A)^{-1}$.
Moreover, $A$ is monotone if 
\begin{equation}
(\forall(x,y)\in\HH\times\HH)
(\forall(u,v)\in Ax\times Ay)\;\;\scal{x-y}{u-v}\geq 0,
\end{equation}
and maximally monotone if it is monotone and there exists no 
monotone operator 
$B\colon\HH\to2^\HH$ such that $\gra A\subset\gra B$ and $A\neq B$.
An operator $B\colon\HH\to\HH$ is $\beta$-cocoercive for some 
$\beta\in\RPP$ if
\begin{equation}
\label{e:cocoercive}
(\forall x\in\HH)(\forall y\in\HH)\quad
\scal{x-y}{Bx-By}\geq\beta\|Bx-By\|^2.
\end{equation}
The conjugate of a function $f\colon\HH\to\RX$ is 
\begin{equation}
\label{e:conjugate}
f^*\colon\HH\to\RXX\colon u\mapsto \sup_{x\in\HH}\big(\scal{x}{u}-f(x)\big), 
\end{equation}
and the infimal convolution of $f$ with $g\colon\HH\to\RX$ is
\begin{equation}
f\infconv g\colon\HH\to\RXX\colon x\mapsto
\inf_{y\in\HH}\big(f(y)+g(x-y)\big).
\end{equation}
The class of lower semicontinuous convex functions $f\colon\HH\to\RX$
such that $\dom f=\menge{x\in\HH}{f(x)<\pinf}\neq\emp$ is denoted by 
$\Gamma_0(\HH)$. If $f\in\Gamma_0(\HH)$, then $f^*\in\Gamma_0(\HH)$
and the subdifferential of $f$ is the maximally monotone operator
\begin{align}
\partial f\colon&\HH\to 2^{\HH}\nonumber\\ 
&x \mapsto\menge{u\in\HH}{(\forall y\in\HH)\;
\scal{y-x}{u}+f(x)\leq f(y)}.
\end{align} 
Let $U\in\BP_{\alpha}(\HH)$ for some $\alpha\in\RPP$. 
The proximity operator of $f\in\Gamma_0(\HH)$ relative to 
the metric induced by $U$ is \cite[Section~XV.4]{Hiriart_Urruty_1996_book_2_convex_amaIf} 
\begin{equation}
\prox^{U}_f\colon\HH\to\HH\colon x
\mapsto\underset{y\in\HH}{\argmin}\:f(y)+\frac12\|x-y\|_{U}^2.
\end{equation}
When $U = \Id$, we retrieve the standard definition of the proximity 
operator \cite{Bauschke_H_2011_book_con_amo,Moreau_J_1965_bsmf_Proximite_eddueh}.
Let $C$ be a nonempty subset of $\HH$. 
The indicator function 
of $C$ 
is defined on $\HH$ as 
\begin{align}
&\iota_C\colon x\mapsto
\begin{cases}
0,&\text{if}\;\;x\in C;\\
\pinf,&\text{if}\;\;x\notin C.
\end{cases}
\end{align}
Finally, $\ell_+^1(\NN)$ denotes the set of summable sequences in 
$\RP$.

\section{A general form of Forward-Backward algorithm}
\label{sec:3}

Optimization problems can often be reduced to finding
a zero of a sum of two maximally monotone operators $A$
and $B$ acting on $\HH$. When $B$ is 
cocoercive (see \eqref{e:cocoercive}), 
a useful algorithm to solve this problem is the forward-backward
algorithm,
which can be formulated in a general form involving a variable
metric as shown in the next result.

\begin{theorem}
\label{t:1}
Let $\alpha\in\RPP$, let $\beta\in\RPP$, let $A\colon\HH\to 2^{\HH}$ be maximally monotone,
and let $B\colon\HH\to\HH$ be cocoercive.
Let $(\eta_n)_{n\in\NN}\in\ell_+^1(\NN)$,
and let $(V_n)_{n\in\NN}$ be a sequence in 
$\BP_{\alpha}(\HH)$ such that
\begin{equation}
\label{e:palawan2012-03-09}
\begin{cases}
\sup_{n\in\NN}\|V_n\|<\pinf\\
(\forall n\in\NN)\quad (1+\eta_n)V_{n+1} \succcurlyeq V_n
\end{cases}
\end{equation}
and $V_n^{1/2} B V_n^{1/2}$ is $\beta$-cocoercive.
Let $(\lambda_n)_{n\in\NN}$ be a sequence in 
$\left]0,1\right]$ such that $\inf_{n\in\NN}\lambda_n >0$
and let $(\gamma_n)_{n\in\NN}$ be a 
sequence in $]0,2\beta[$ such that $\inf_{n\in\NN} \gamma_n > 0$
and $\sup_{n\in \NN} \gamma_n < 2\beta$.
Let $x_0\in\HH$, and let $(a_n)_{n\in\NN}$ 
and $(b_n)_{n\in\NN}$ be 
absolutely summable sequences in $\HH$.
Suppose that
$Z=\zer(A+B)\neq\emp$,
and set 
\begin{equation}
\label{e:forward}
(\forall n\in\NN)\quad
\begin{array}{l}
\left\lfloor
\begin{array}{l}
y_n= x_n-\gamma_n V_n (Bx_n+b_n)\\[1mm]
x_{n+1}=x_n+\lambda_n\big(J_{\gamma_n V_n A}\,
(y_n)+a_n-x_n\big).
\end{array}
\right.\\[2mm]
\end{array}
\end{equation}
Then $x_n\weakly\overline{x}$ for some $\overline{x}\in Z$.
\end{theorem}
At iteration $n$, variables $a_n$ and $b_n$ model 
numerical errors possibly arising when applying $J_{\gamma_n V_n
A}$ or $B$. Note also that, if $B$ is $\mu$-cocoercive with $\mu \in \RPP$,
one can choose $\beta = \mu (\sup_{n\in \NN}\|V_n\|)^{-1}$, 
which allows us to retrieve \cite[Theorem~4.1]{Combettes_P_2014_j-optim_Variable_mfb}.
In the next section, we shall see how a judicious use of
this result allows us to derive a variety of flexible convex
optimization algorithms.  

\section{A variable metric primal-dual method}
\label{sec:4}
\subsection{Formulation}
A wide array of optimization problems encountered in image processing
are instances of the following one, which was first investigated in
\cite{Combettes_P_2012_j-svva_pri_dsa} and can be viewed as a more
structured version of the minimization problem in \eqref{e:prob1}:
\begin{problem}
\label{prob:2}
Let  $z\in\HH$, let $m$ be a 
strictly positive integer, let $f\in\Gamma_0(\HH)$, and let
$h\colon\HH\to\RR$ be convex and differentiable with a 
Lipschitzian gradient.
For every $i\in\{1,\ldots,m\}$, let $r_i\in\GG_i$, let 
$g_i\in\Gamma_0(\GG_i)$, let $\ell_i\in\Gamma_0(\GG_i)$ 
be strongly convex,\footnote{For every $i\in\{1,\ldots,m\}$,
$\ell_i$ is $\nu_i^{-1}$-strongly convex with $\nu_i\in \RPP$
if and only if $\ell_i^*$ is $\nu_i$-Lipschitz differentiable
\cite[Theorem~18.15]{Bauschke_H_2011_book_con_amo}.} 
and suppose that $0\neq L_i\in\BL(\HH,\GG_i)$. 
Suppose that
\begin{equation}
\label{e:2012-30-04}
z\in\ran\bigg(\partial f+\sum_{i=1}^mL_i^*(\partial g_i\infconv
\partial\ell_i)(L_i\cdot-r_i)+\nabla h\bigg).
\end{equation}
Consider the problem
\begin{equation}
\label{e:primalvar}
\minimize{x\in\HH}{f(x)+\sum_{i=1}^m\,(g_i\infconv\ell_i)
(L_ix-r_i)+h(x)-\scal{x}{z}}, 
\end{equation}
and the dual problem
\begin{align}
\label{e:dualvar}
\minimize{v_1\in\GG_1,\ldots,v_m\in\GG_m}{}
&\big(f^*\infconv h^*\big)\bigg(z-\sum_{i=1}^mL_i^*v_i\bigg)\nonumber\\
&+\sum_{i=1}^m\big(g_i^*(v_i)+\ell_i^*(v_i)+\scal{v_i}{r_i}\big).
\end{align}
\end{problem}

Note that in the special case when $\ell_i = \iota_{\{0\}}$, 
$g_i\infconv\ell_i$ reduces to $g_i$ in \eqref{e:primalvar}.

Let us now examine how Problem \ref{prob:2} can be reformulated from the standpoint of monotone operators.
To this end, let us define $\boldsymbol{g} \in \Gamma_0(\GGG)$, $\boldsymbol{\ell}\in \Gamma_0(\GGG)$
and $\boldsymbol{L}\in \BL(\HH,\GGG)$ by
\begin{align}
&\boldsymbol{g}\colon\boldsymbol{v} \mapsto \sum_{i=1}^m g_i(v_i), 
\qquad \boldsymbol{\ell}\colon\boldsymbol{v} \mapsto \sum_{i=1}^m \ell_i(v_i)\nonumber\\
\text{and}\quad
&\boldsymbol{L}\colon x \mapsto (L_1 x,\ldots,L_m x).
\end{align}
Let us now introduce the product space $\KKK=\HH\oplus\GGG$ and 
the operators
\begin{align}
\label{e:maximal1}
{\boldsymbol A}\colon\quad\;\;\KKK&\to 2^{\KKK}\nonumber\\
(x,\boldsymbol{v})&\mapsto (\partial f(x)-z+\boldsymbol{L}^*
\boldsymbol{v}) \times(-\boldsymbol{L} x+\partial
\boldsymbol{g}^*(\boldsymbol{v})+\boldsymbol{r})
\end{align}
and
\begin{align}
\label{e:maximal21}
{\boldsymbol B}\colon\quad\;\;\KKK&\to\KKK\nonumber\\
(x,\boldsymbol{v})&\mapsto 
\big(\nabla h(x),\nabla \boldsymbol{\ell}^*(\boldsymbol{v})\big).
\end{align}
The operator ${\boldsymbol A}$ can be shown to be maximally
monotone,whereas ${\boldsymbol B}$ is cocoercive.
A key observation in this context is that, if there exists
$(\overline{x},\overline{\boldsymbol{v}})
\in \KKK$ such that $(\overline{x},\overline{\boldsymbol{v}}) \in
\zer({\boldsymbol A}+{\boldsymbol B})$,
then $(\overline{x},\overline{\boldsymbol{v}})$ is a pair of
primal-dual solutions to
Problem \ref{prob:2} \cite{Combettes_P_2012_j-svva_pri_dsa}. This
connection with the construction for a zero of ${\boldsymbol
A}+{\boldsymbol B}$
makes it possible to apply a forward-backward algorithm as
discussed in Section \ref{sec:3}, by using
a linear operator $\boldsymbol{V}_n \in \BL(\KKK,\KKK)$ to change
the metric at each iteration $n$.
Depending on the form of this operator various algorithms can be obtained.

\subsection{A first class of primal-dual algorithms}
Let $\alpha\in\RPP$, let $(U_n)_{n\in\NN}$ be a sequence 
in $\BP_{\alpha}(\HH)$ such that 
$(\forall n\in\NN)$ $U_{n+1}\succcurlyeq U_n$.
For every $i\in\{1,\ldots,m\}$, 
let  $(U_{i,n})_{n\in\NN}$ be a sequence in $\BP_{\alpha}(\GG_i)$
such that $(\forall n\in\NN)$ $U_{i,n+1}\succcurlyeq U_{i,n}$.
A first possible choice for $(\boldsymbol{V}_n)_{n\in\NN}$ is given by
\begin{align}
(\forall n \in \NN)\quad
\boldsymbol{V}_n^{-1}\colon 
(x,\boldsymbol{v}) &\mapsto (U_n^{-1} x-\boldsymbol{L}^* \boldsymbol{v},-\boldsymbol{L}x+\widetilde{\boldsymbol{U}}_n^{-1}\boldsymbol{v})
\end{align}
where
\begin{equation}\label{e:defUt}
\widetilde{\boldsymbol{U}}_n\colon \GGG\to \GGG\colon (v_1,\ldots,v_m) \mapsto (U_{1,n} v_1,\ldots,U_{m,n} v_m).
\end{equation}
The following result constitutes a direct extension of \cite[Example~6.4]{Combettes_P_2014_j-optim_Variable_mfb}:
\begin{proposition}
\label{ex:2012-30-04}
Let $x_0\in \HH$, and let $(a_n)_{n\in\NN}$ and $(c_n)_{n\in\NN}$ be absolutely summable 
sequences in $\HH$.
For every $i\in\{1,\ldots,m\}$, 
let $v_{i,0}\in\GG_i$, let $(b_{i,n})_{n\in\NN}$ and
$(d_{i,n})_{n\in\NN}$ be absolutely summable sequences in $\GG_i$.
For every $n\in \NN$, 
let  $\mu_n \in \RPP$ be a Lipschitz constant of $U_n^{1/2}\circ \nabla h \circ U_n^{1/2}$
and, for every $i\in \{1,\ldots,m\}$, let $\nu_{i,n}\in \RPP$ be a Lipschitz
constant of $U_{i,n}^{1/2} \circ \nabla \ell_i^*\circ  U_{i,n}^{1/2}$.
Let $(\lambda_n)_{n\in\NN}$ be a sequence in $]0,1]$
such that $\inf_{n\in \NN} \lambda_n > 0$.
For every $n\in\NN$, set
\begin{equation}
\label{e:121}
\delta_n=
\Bigg(\sum_{i= 1}^m\|\sqrt{U_{i,n}} 
L_i\sqrt{U_n}\|^2\Bigg)^{-1/2}-1,
\end{equation}
and suppose that
\begin{equation}
\label{e:2f9h79p}
\inf_{n\in \NN} \frac{\delta_n}
{(1+\delta_n)\max\{\mu_n,\nu_{1,n},\ldots,\nu_{m,n}\}}
>
\frac{1}{2}.
\end{equation}
Set
\begin{align}
\label{e:cocoeqal8:20coco}
&\operatorname{For}\;n=0,1,\ldots\nonumber\\
&\begin{array}{l}
\left\lfloor
\begin{array}{l}
p_n=\prox^{U^{-1}_n}_f\Big(x_n-U_n
\big(\sum_{i=1}^{m}L_{i}^*v_{i,n}+
\nabla h(x_n)\\
\qquad\qquad\qquad\;\;+c_n-z\big)\Big)+a_n\\
y_n=2p_n-x_n\\
x_{n+1}=x_n+\lambda_n(p_n-x_n)\\
\operatorname{For}\;  i=1,\ldots, m\\
\left\lfloor
\begin{array}{l}
q_{i,n}=\prox^{U_{i,n}^{-1}}_{g_{i}^{*}}
\Big(v_{i,n}+U_{i,n}\big(L_iy_n- 
\nabla\ell_{i}^{*}(v_{i,n})\\
\qquad\qquad\qquad\;\;\;\;-d_{i,n}-r_i\big)\Big)+b_{i,n}\\
v_{i,n+1}=v_{i,n}+\lambda_n(q_{i,n}-v_{i,n}).\\
\end{array}
\right.\\[2mm]
\end{array}
\right.\\[2mm]
\end{array}
\end{align}
Then $(x_n)_{n\in\NN}$ converges weakly to a solution to 
\eqref{e:primalvar}, for every $i\in\{1,\ldots,m\}$ 
$(v_{i,n})_{n\in\NN}$ converges weakly to some
$\overline{v}_i\in\GG_i$, and 
$(\overline{v}_{1},\ldots,\overline{v}_{m})$ is 
a solution to \eqref{e:dualvar}.
\end{proposition}

In the special case when $U_n \equiv \tau \Id$ with $\tau \in \RPP$ and, for every $i\in \{1,\ldots,m\}$,
$U_{i,n} \equiv \sigma_i \Id$ with $\sigma_i \in \RPP$, we recover the parallel algorithm proposed in \cite{Vu_B_2013_j-acm_spl_adm}.
Variants of this algorithm where, for every $i\in \{1,\ldots,m\}$, $\ell_i = \iota_{\{0\}}$ are
also investigated in \cite{Condat_L_2013_j-ota-primal-dsm}. In this case, less restrictive assumptions on the choice
of $(\tau,\sigma_1,\ldots,\sigma_m)$ can be made. Note that this algorithm itself can be viewed as a generalization
of the algorithm which constitutes the main topic of \cite{Chambolle_A_2010_first_opdacpai,Esser_E_2010_j-siam-is_gen_fcf,He_B_2012_j-siam-is_conv_apd} (designated by some authors as PDHG). 
A preconditioned version of this algorithm was proposed in \cite{Pock_T_2008_p-iccv_diagonal_pffo} 
corresponding to the case when $m=1$,
$(\forall n\in\NN)$ $U_n$ and $U_{1,n}$ are constant matrices, and no error term is taken into account.
Algorithm \eqref{e:cocoeqal8:20coco} when, for every $n\in \NN$, $\lambda_n \equiv 1$, $U_n$ and $(U_{i,n})_{1\le i \le m}$ are diagonal matrices,
$h = 0$, and $(\forall i \in \{1,\ldots,m\})$ $\ell_i = \iota_{\{0\}}$ appears also to be closely related
to the adaptive method in \cite{Goldstein_T_2013_adaptive_pdh}.

\subsection{A second class of primal-dual algorithms}
Let $\alpha\in\RPP$, let $(U_n)_{n\in\NN}$ be a sequence 
in $\BP_{\alpha}(\HH)$ such that 
$(\forall n\in\NN)$ $U_{n+1}\succcurlyeq U_n$.
For every $i\in\{1,\ldots,m\}$, 
let  $(U_{i,n})_{n\in\NN}$ be a sequence in $\BP_{\alpha}(\GG_i)$
such that $(\forall n\in\NN)$ $U_{i,n+1}\succcurlyeq U_{i,n}$.
A second possible choice for $(\boldsymbol{V}_n)_{n\in\NN}$ is given by
the following diagonal form:
\begin{align}
(\forall n \in \NN)\quad
\boldsymbol{V}_n^{-1}\colon 
(x,\boldsymbol{v}) &\mapsto \big(U_n^{-1} x,(\widetilde{\boldsymbol{U}}_n^{-1} - \boldsymbol{L}U_n \boldsymbol{L}^*)\boldsymbol{v} \big)
\end{align}
where $\widetilde{\boldsymbol{U}}_n$ is given by \eqref{e:defUt}.

The following result can then be deduced from Theorem~\ref{t:1}.
Its proof is skipped due to the lack of space.

\begin{proposition}
\label{ex:2014-25-01}
Let $x_0\in \HH$, and let $(c_n)_{n\in\NN}$ be an absolutely summable 
sequence in $\HH$.
For every $i\in\{1,\ldots,m\}$, let $v_{i,0} \in \GG_i$,
let $(b_{i,n})_{n\in\NN}$ and
$(d_{i,n})_{n\in\NN}$ be absolutely summable sequences in $\GG_i$.
For every $n\in \NN$, 
let  $\mu_n \in \RPP$ be a Lipschitz constant of $U_n^{1/2}\circ \nabla h \circ U_n^{1/2}$
and, for every $i\in \{1,\ldots,m\}$, let $\nu_{i,n}\in \RPP$ be a Lipschitz
constant of $U_{i,n}^{1/2} \circ \nabla \ell_i^*\circ  U_{i,n}^{1/2}$.
Let $(\lambda_n)_{n\in\NN}$ be a sequence in 
$]0,1]$ such that $\inf_{n\in \NN} \lambda_n > 0$.
For every $n\in\NN$, set
\begin{equation}
\zeta_n= 1-\sum_{i= 1}^m\|\sqrt{U_{i,n}} L_i\sqrt{U_n}\|^2
\end{equation}
and suppose that 
\begin{equation}
\inf_{n\in \NN} \frac{\zeta_n}
{\max\{\zeta_n\mu_n,\nu_{1,n},\ldots,\nu_{m,n}\}}
>\frac{1}{2}.
\end{equation}
Set
\begin{align}
\label{e:zhangcoco}
&\operatorname{For}\;n=0,1,\ldots\nonumber\\
&\begin{array}{l}
\left\lfloor
\begin{array}{l}
s_n=x_n - U_n (\nabla h(x_n)+c_n-z)\\
y_n = s_n - U_n \sum_{i=1}^m L_i^* v_{i,n}\\
\operatorname{For}\;  i=1,\ldots, m\\
\left\lfloor
\begin{array}{l}
q_{i,n}=\prox^{U_{i,n}^{-1}}_{g_{i}^{*}}
\Big(v_{i,n}+U_{i,n}\big(L_iy_n- 
\nabla\ell_{i}^{*}(v_{i,n})\\
\qquad\qquad\qquad\;\;\;\;-d_{i,n}-r_i\big)\Big)+b_{i,n}\\
v_{i,n+1}=v_{i,n}+\lambda_n(q_{i,n}-v_{i,n}).\\
\end{array}
\right.\\[2mm]
p_n = s_n-U_n \sum_{i=1}^m L_i^* q_{i,n}\\
x_{n+1}=x_n+\lambda_n(p_n-x_n).\\
\end{array}
\right.\\[2mm]
\end{array}
\end{align}
Assume that $f = 0$.
Then $(x_n)_{n\in\NN}$ converges weakly to a solution to 
\eqref{e:primalvar}, for every $i\in\{1,\ldots,m\}$ 
$(v_{i,n})_{n\in\NN}$ converges weakly to some
$\overline{v}_i\in\GG_i$, and 
$(\overline{v}_{1},\ldots,\overline{v}_{m})$ is 
a solution to \eqref{e:dualvar}.
\end{proposition}
The algorithm proposed in \cite{Loris_I_2011_generalization_ist,Chen_P_2013_j-inv-prob_prim_dfp}
is a special case of the previous one, in the absence of errors, when $m = 1$, 
$\HH$ and $\GG_1$ are finite dimensional spaces,
$\ell_1 = \iota_{\{0\}}$,
$U_n \equiv \tau \Id$ with $\tau \in \RPP$,
$U_{1,n} \equiv \sigma \Id$ with $\sigma \in \RPP$, and
no relaxation ($\lambda_n \equiv 1)$ or a constant one ($\lambda_n \equiv \kappa < 1$) is performed.

\begin{figure}[t]
\begin{center}
\begin{tabular}{@{}c@{}c@{}}
\includegraphics[width=4.3cm]{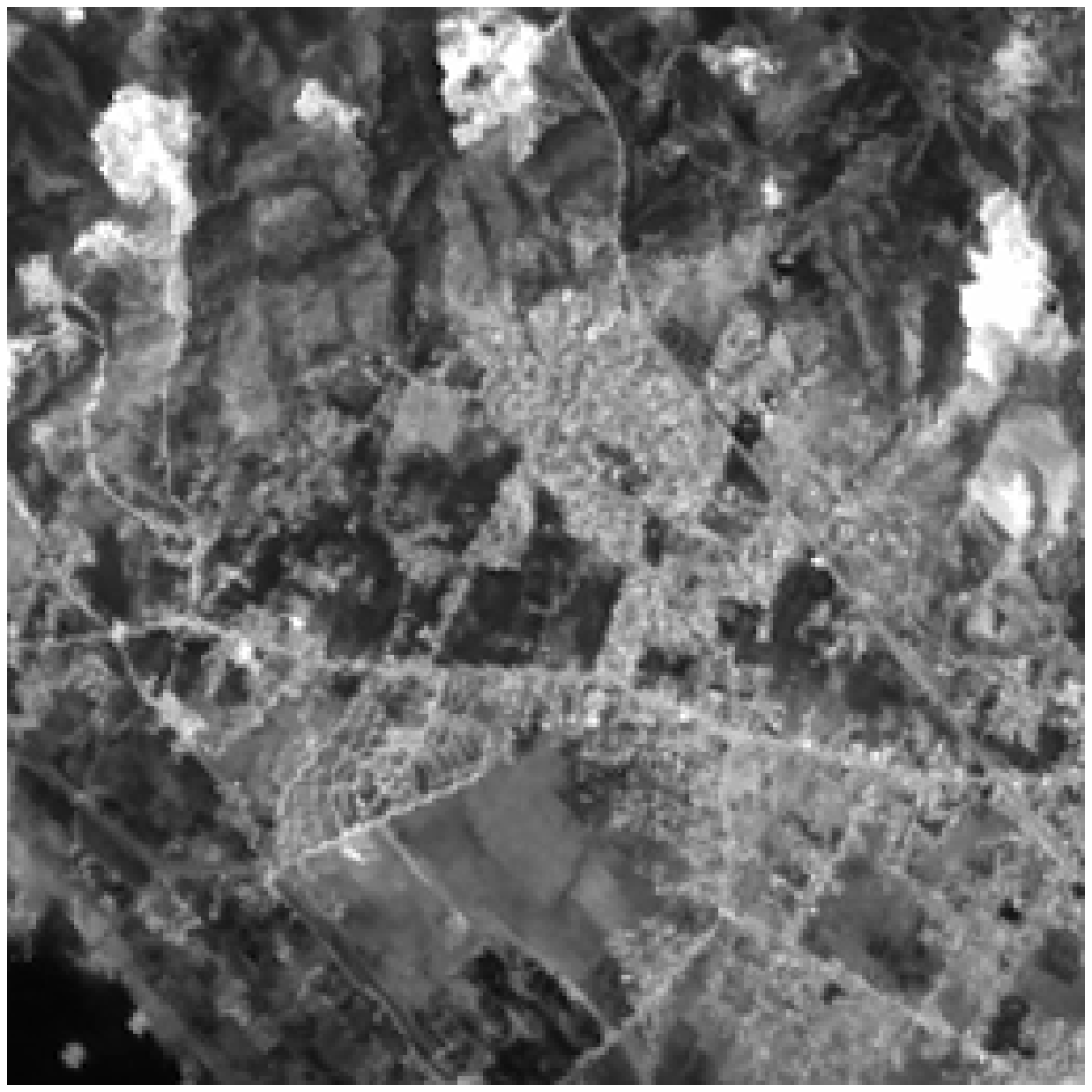}&
\includegraphics[width=4.3cm]{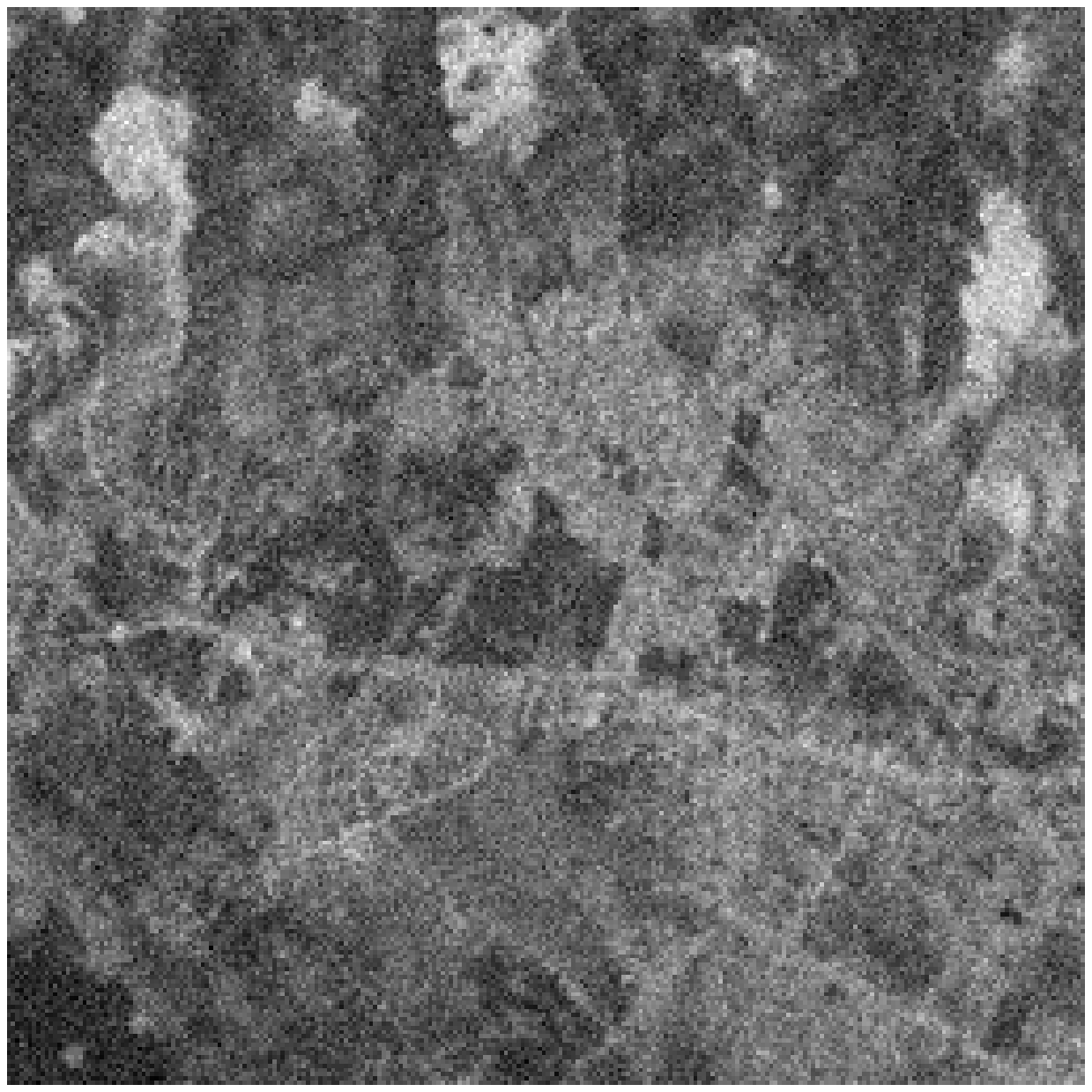}\\
\small{(a)} & \small{(b)}\\
\includegraphics[width=4.3cm]{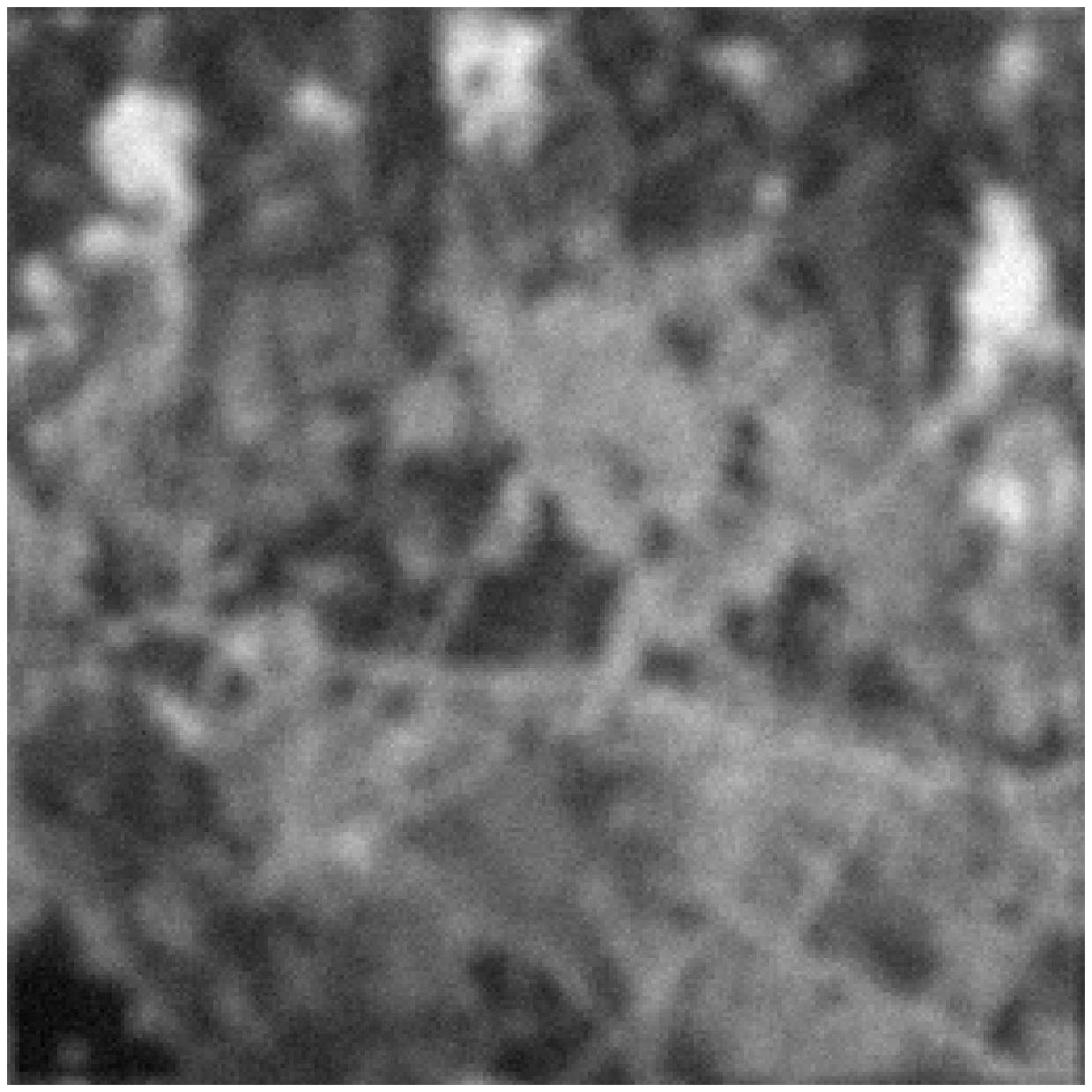}&
\includegraphics[width=4.3cm]{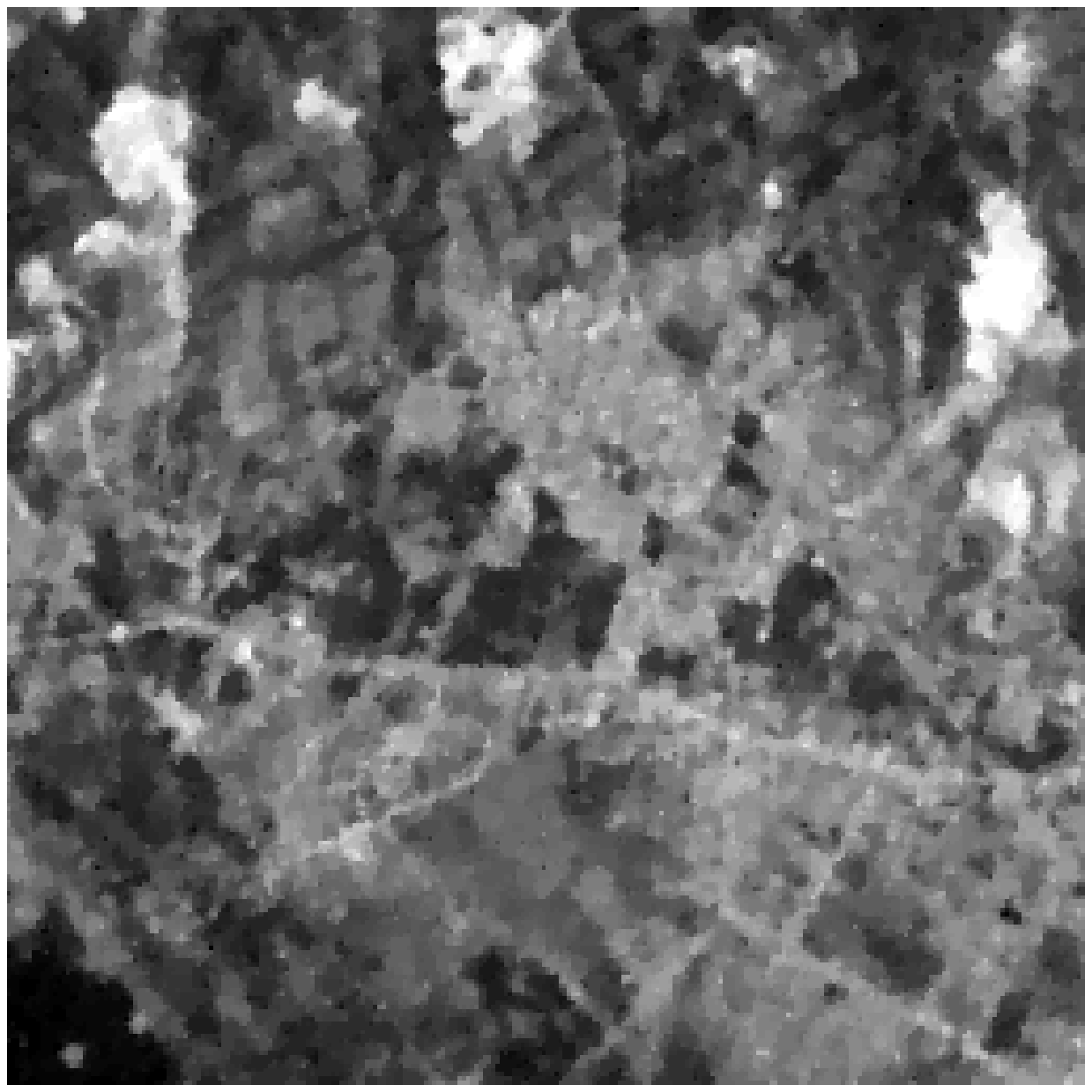}\\
\small{(c)} & \small{(d)}
\end{tabular} 
\caption{Original image $\overline{x}$ (a), noisy image $w_1$ (SNR = $5.87$~dB) (b),
blurred image $w_2$ (SNR = $16.63$~dB) (c), and restored image
$\widetilde{x}$  (SNR = $21.61$~dB)  (d).}
\label{fig:rest}
\end{center}
\end{figure}

\section{Application to image restoration}
\label{sec:5}
We illustrate the flexibility of the proposed primal-dual algorithms
on an image recovery example. Two observed images $w_1$ and $w_2$ of the same scene
$\overline{x}\in \RR^{N}$ ($N=256^2$)
are available (see Fig.~\ref{fig:rest}(a)-(c)). The first one is
corrupted with a noise with a variance
$\theta_1^2 = 576$, while the second one has been degraded by a linear operator $H\in \RR^{N\times N}$
($7\times 7$ uniform blur) and a noise with variance $\theta_2^2 = 25$. The noise components
are mutually statistically independent, additive, zero-mean, white, and Gaussian distributed.
Note that this kind of multivariate restoration problem is
encountered in some push-broom satellite imaging systems.

An estimate $\widetilde{x}$ of $\overline{x}$ is computed as a solution to \eqref{e:primalvar}
where $m=2$, $z=0$, $r_1=0$, $r_2 = 0$,
\begin{align}
h &= \frac{1}{\theta_1^2} \|\cdot-w_1\|^2+ \frac{1}{\theta_2^2}\|H\cdot-w_2\|^2,\\
g_1 & = \iota_{[0,255]^N},\quad  g_2 = \kappa \|\cdot\|_{1,2}, \label{e:g2}\\
f &= 0,\quad \ell_1 = \ell_2 = \iota_{\{0\}}
\end{align}
where the second function in \eqref{e:g2} denotes the $\ell_{1,2}$-norm and $\kappa \in \RPP$.
In addition, $L_1 = \Id$ and $L_2 = [G_1^\top,G_2^\top]^\top$ where $G_1\in \RR^{N\times N}$ and $G_2^{N\times N}$ are
horizontal and vertical discrete gradient operators. Function $g_1$ introduces some a priori
constraint on the range values in the target image, while function $g_2\circ L_2$ corresponds to
a classical total variation regularization. The minimization problem
is solved numerically by using Algorithm \eqref{e:zhangcoco} with $\lambda_n \equiv 1$. 
In a first experiment,
standard choices of the algorithm parameters are made by setting $U_n \equiv \tau \Id$,
$U_{1,n} \equiv \sigma_1 \Id$, and $U_{2,n} = \sigma_2 \Id$ with $(\tau,\sigma_1,\sigma_2)
\in \RPP^3$. In a second experiment, a more sophisticated choice of the metric is made.
The operators $(U_n)_{n\in \NN}$, $(U_{1,n})_{n\in \NN}$ and $(U_{2,n})_{n\in \NN}$ are still
chosen diagonal and constant in order to facilitate the implementation of the algorithm,
but the diagonal values are optimized in an empirical manner. A similar strategy was applied
in \cite{Pock_T_2008_p-iccv_diagonal_pffo} in the case of Algorithm \eqref{e:cocoeqal8:20coco}. 
The regularization parameter $\kappa$ has been
set so as to get the highest value of the resulting signal-to-noise ratio (SNR).

\begin{figure}
\centering
\includegraphics[width=7.5cm]{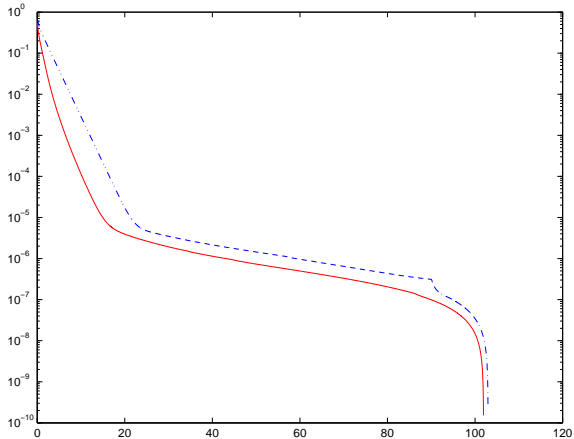}
\caption{Normalized norm of the error on the iterate vs computation time (in seconds)
for Experiment 1 (blue, dash dot line) and Experiment 2 (red, continuous line).}
\label{fig:conv}
\end{figure}

The restored image is displayed in Fig.~\ref{fig:rest}(d).
Fig.~\ref{fig:conv} shows the convergence profile of the algorithm.
We plot the evolution of the normalized Euclidean distance (in log scale) 
between the iterates  and $\widetilde{x}$ in terms of computational time 
(\texttt{Matlab R2011b} 
codes running on a single-core Intel i7-2620M CPU@2.7 GHz with 8 GB of RAM).
An approximation of $\widetilde{x}$ obtained after 5000 iterations is used.
This result illustrates the fact that an appropriate choice of the
metric may be beneficial in terms of speed of convergence.

\end{document}